\documentclass [11pt,a4paper]{amsart}
\usepackage[english]{babel}
\newcounter{lemma}
\newtheorem{thm}{Theorem}

\newtheorem{lem}[lemma]{Lemma}

\newcommand{\de}{\partial}

\newcommand{\R}{\mathbb{R}}

\newcommand{\Ca}{\mathrm{Ca}}

\newcommand{\ve}{\varepsilon}

\newcommand{\jb}{\bar\jmath}

\newcommand{\Ric}{\mathrm{Ric}}
\newcommand{\ddbar}{\frac{\sqrt{-1}}{2} \partial \overline{\partial}}
\newcommand{\ddt}{\frac{\partial}{\partial t}}

\newenvironment{proof1}{\noindent\textbf{Proof.} }{\hfill
\medskip}

\begin{document}

\title[THE CALABI FLOW]{THE CALABI FLOW WITH SMALL INITIAL ENERGY}
\maketitle

\begin{center}
\begin{tabular}{ccc}
{\bf Valentino Tosatti}  & & {\bf Ben Weinkove} \footnotemark \\
Harvard University & & Harvard University \\
Department of Mathematics & & Department of Mathematics \\
Cambridge, MA 02138 & & Cambridge, MA 02138
\end{tabular}
\end{center}
\footnotetext{Most of this work was carried out while the second author was visiting Imperial College on a  Royal
Society Research Assistantship;  the second author is also supported in part by
National Science Foundation grant DMS-05-04285}
\bigskip
\bigskip

{\bf Abstract.}  We show that on K\"ahler manifolds $M$ with $c_1(M)=0$
 the Calabi flow converges to a constant scalar
curvature metric if the initial Calabi energy is sufficiently small.
We prove a similar result on manifolds with $c_1(M)<0$ if the
K\"ahler class is close to the canonical class.

\bigskip

\section{Introduction}

\bigskip

It has been known since  Yau's proof of the Calabi conjecture
 that there always exists a K\"ahler-Einstein metric on a compact K\"ahler manifold $M$ with
 $c_1(M)=0$ \cite{yau1}
and $c_1(M)<0$ \cite{yau1}, \cite{aubin}.  The problem is still open
in the remaining case $c_1(M)>0$, where existence was conjectured by
Yau \cite{yau2} to be equivalent to stability in the sense of
geometric invariant theory. There has been much recent progress
towards understanding and refining this conjecture and it is now
expected, more generally, that the existence of a K\"ahler metric of
constant scalar curvature should be equivalent to some notion of
K-stability \cite{Ti}, \cite{donaldson2002}.
\\

A natural approach to the problem of existence of these canonical metrics
is to consider the Calabi flow.  Fix a K\"ahler metric $\omega = \frac{\sqrt{-1}}{2} g_{i \overline{j}} dz^i \wedge dz^{\overline{j}}$.  Then the Calabi flow is
a flow of K\"ahler potentials $\varphi
=\varphi(t)$
given by
$$\ddt{} \varphi = R(\omega_{\varphi}) - \mu,$$
where $\omega_{\varphi} = \omega + \ddbar \varphi>0$,  $R(\omega_{\varphi})$
is the scalar curvature of the metric $\omega_{\varphi}$, and
 $$\mu  = \frac{\pi [c_1(M)] \cdot [\omega]^{n-1}}{[\omega]^n}, \quad  n = \textrm{dim}_{\mathbb{C}} M$$
It is a gradient flow for the Mabuchi energy \cite{mabuchi}, and it decreases the Calabi
functional, which we define by
$$\Ca(\omega_{\varphi})=\int_M(R(\omega_{\varphi})-\mu)^2\omega_{\varphi}^n,$$
(this differs from the usual definition by a constant). It is
expected that the Calabi flow should converge to a constant scalar
curvature (cscK) metric when one exists \cite{donaldson}.
Unfortunately, since the flow is a parabolic fourth order equation,
there are very few techniques available to tackle this kind of
analytic question. In the case of Riemann surfaces, Chru\'sciel
\cite{chrusciel} showed that the flow always converges to a constant
curvature metric, making use of the existence of such a metric and
the Bondi mass of general relativity. Chen \cite{chen2} (see also
\cite{struwe}) gave a different proof of Chru\'sciel's theorem using
concentration-compactness results.  The analytic difficulties become
more acute in higher dimensions, but some progress has been made
\cite{chang}, \cite{che}.  In particular, Chen and He \cite{che}
recently showed that the Calabi flow exists as long as the Ricci
curvature remains bounded, and proved a short time existence result.
They also showed that if the initial metric is close to a cscK
metric in the $C^{3+{\alpha}}$ sense, for some $\alpha>0$,
then the flow will converge to the cscK metric. \\

In this paper we restrict to  manifolds satisfying $c_1(M)=0$ and
$c_1(M)<0$ and consider the case when the initial data has small
Calabi energy, in a suitable sense.  When $c_1(M)=0$ every class
admits a unique K\"ahler Ricci-flat metric by Yau's theorem.  We have the
following result.

\begin{thm}\label{main1}
Let $M$ be a compact K\"ahler manifold with $c_1(M)=0$. Fix constants $K_1,K_2$ and a reference K\"ahler metric
$\omega$.
Then there exists $\ve>0$ depending only on $\omega, K_1, K_2$  such that
if $\omega_0$ is any  K\"ahler metric in $[\omega]$ satisfying:
\begin{enumerate}
\item[(i)] $\displaystyle{-K_1 \omega \le \Ric(\omega_0) \le K_2 \omega}$
\item[(ii)] $\displaystyle{ \Ca(\omega_0) < \ve}$,
\end{enumerate}
then the solution of the Calabi flow starting at $\omega_0$ exists for all time
and converges in $C^{\infty}$ to the K\"ahler Ricci-flat metric in $[\omega]$.
\end{thm}

In the case $c_1(M)<0$ we know that there exists a K\"ahler-Einstein
metric in the class $-\pi c_1(M)$, and, by the implicit function
theorem, there exist cscK metrics in nearby classes \cite{lesim}.  Indeed, it is
expected by the results of \cite{chen0}, \cite{we}, \cite{sowe} that
there exist cscK metrics in a larger neighborhood of $-\pi
c_1(M)$ (but not necessarily in every class \cite{ross}).
  We show that the Calabi flow converges to the unique cscK metric if the
  Calabi energy is small and if the class is sufficiently close to $- \pi
  c_1(M)$.\\

  Fix a norm $\| \cdot \|_N$ on the vector space $H^{1,1}_\R(M)$.
  Then we have the following:

\begin{thm}\label{main2}
Let $M$ be a compact K\"ahler manifold with $c_1(M)<0$. Fix two constants $K_1,K_2$ and a reference metric
$\omega \in - \pi c_1(M)$ with $\Ric(\omega)<0$.
Then there exist constants $\ve, \delta >0$ depending only on $\omega, K_1, K_2$ and $\|
\cdot \|_N$ such that
if a K\"ahler metric $\omega_0$ satisfies the conditions:

\begin{enumerate}
\item[(i)] $\displaystyle{\| [\omega_0] - [\omega] \|_N < \delta}$
\item[(ii)] $\displaystyle{R (\omega_0) \ge - K_1}$
\item[(iii)] $\displaystyle{\Ric(\omega_0) \le K_2 \omega}$
\item[(iv)] $\displaystyle{ \Ca(\omega_0) < \ve}$,
\end{enumerate}
then the solution of the Calabi flow starting at $\omega_0$ exists for all time
and converges in $C^{\infty}$ to a cscK metric in $[\omega_0]$.
\end{thm}

Note that we do not assume \emph{a priori} that a cscK metric exists
in the given class.  The results here are quite different from the stability result of \cite{che} in these cases, where it is
assumed that a cscK metric exists and that the initial metric is
close to the cscK metric in the $C^{3+\alpha}$ norm.\\

In section 2, we prove two key compactness lemmas for K\"ahler metrics, which rely on Yau's estimates.
These, along with the `short-time existence estimates' of Chen-He and the fact that the Calabi functional is decreasing,  are used in section 3 to prove Theorems 1 and 2.  In the course of the proofs, we will use the letter $C$ to denote a suitably uniform positive constant, which may differ from line to line.

\section{Compactness results}

In this section we prove two compactness results for K\"ahler
metrics.  The first is as follows.

\begin{lem} \label{compactnessthm1}
Let $M$ be a compact K\"ahler manifold.  Fix a reference K\"ahler
metric $\omega$ and suppose $\omega' = \omega + \ddbar \varphi$ is a
metric satisfying
$$-K_1 \omega \le \Ric(\omega') \le K_2 \omega.$$
Then for any $0< \alpha<1$ there exist positive constants $c$ and $A$  depending only on $\omega$, $\alpha$, $K_1$ and $K_2$
such that $\omega' \ge c \, \omega$ and $$\| \varphi - \underline{\varphi} \|_{C^{3+\alpha}(\omega)}
\le A,$$
where $\underline{\varphi} = \frac{1}{V} \int_M \varphi \, \omega^n$, for $V= \int_M \omega^n$.
\end{lem}
\begin{proof1}
Define
\begin{equation}\label{monge}
F:=\log\frac{\omega'^n}{\omega^n}.
\end{equation}

Calculate
\begin{equation}\label{lap1}
\Delta F=R-g^{i\jb}R'_{i\jb},
\end{equation}
where we are using the standard notation $R'_{i \jb}$ for the Ricci
curvature of the metric $\omega'$. \\

Our assumptions imply $|\Delta F|\leq C$.
Applying Green's formula for the metric $\omega$ we see that for all
$x \in M$
\begin{equation} \label{greens}
F(x)=\frac{1}{V}\int_M F\omega^n-\int_{y \in M}\Delta
F(y)G(x,y)\omega^n(y),
\end{equation}
where the Green's function $G$ is bounded in $L^1(\omega)$ and
uniformly from below.  Hence
$$ \sup_M |F| \le \left| \frac{1}{V} \int_M F \omega^n \right| + C.$$
It remains then to bound the integral of $F$.  But from the
definition of $F$ there must exist some point $x_0 \in M$ with
$F(x_0)=0$ and so applying (\ref{greens}) with $x=x_0$ shows that
$\int_M F \omega^n$ is bounded above and below.  (In fact, by
Jensen's inequality and the concavity of the log function it is easy
to see that we can take the upper bound to be zero.) \\

Now that we have $|F|\leq C$ and $\Delta F\geq -C$,
Yau's second order estimates \cite[Prop. 2.1]{yau1} give us bounds
$$\sup_M |\varphi-\underline{\varphi}|  \leq C,\hspace{10pt} C^{-1}\omega\leq\omega'\leq C\omega.$$
We also notice that $|\Delta F|\leq C$
together with the elliptic estimates for $\Delta$ give uniform
$L_2^p(\omega)$ bounds on $F$ for all $p\geq 1$. We can now use the
method of Evans and Krylov (see also the simplification by Trudinger
\cite{tr}, and the exposition of Siu  \cite{siu}) that employs
Moser's Harnack inequality to obtain the H\"older estimate $\|
\varphi-\underline{\varphi}\|_{C^{2+\beta}(\omega)}\leq C$ for some
uniform $\beta>0$.  Let $\psi$ be a local potential for $\omega$ so that $\omega' = \ddbar{(\varphi+ \psi)}$ locally.
Differentiate the equation
\eqref{monge} with respect to $z^i$ to obtain
$$\frac{\de F}{\de z^i}=\Delta'\frac{\de (\varphi + \psi) }{\de z^i} -\frac{\de}{\de z^i}\log\det g.$$
The operator $\Delta'$ is uniformly elliptic with coefficients bounded in $C^{\beta}(\omega)$
and $F$ has uniform $L^p_2(\omega)$ bounds for all $p\geq 1$,
so the elliptic estimates give uniform $L^p_4(\omega)$ bounds on $\varphi-\underline{\varphi}$ and Sobolev's embedding theorem
gives uniform $C^{3+\alpha}(\omega)$ bounds after taking $p$ sufficiently large.  Q.E.D.
\end{proof1}

For the case $c_1(M)<0$ we will use instead the following compactness result.

\begin{lem} \label{compactnessthm2}
Let $M$ be a compact K\"ahler manifold with $c_1(M)<0$.  Fix a
reference K\"ahler metric $\omega$ with $\Ric(\omega)<0$ and suppose
$\omega' = \omega +\ddbar \varphi$ is a metric satisfying
$$R(\omega') \ge -K_1 \quad \textrm{and} \quad  \Ric(\omega') \le K_2 \omega.$$
Then for any $0< \alpha<1$ there exist positive constants $c$ and $A$ depending only on $\omega$, $\alpha$, $K_1$ and $K_2$
such that $\omega' \ge c \, \omega$ and $$\| \varphi - \underline{\varphi} \|_{C^{3+\alpha}(\omega)}
\le A.$$
\end{lem}
\begin{proof1}
Using the same notation as in the proof of Lemma
\ref{compactnessthm1} we see immediately that $\Delta F \ge -C$. It
then follows from (\ref{greens}) that $\sup_M F \le C$. Calculate
\begin{equation} \label{Delta'F}
\Delta' F = g'^{i \jb} R_{i \jb} - R(\omega').
\end{equation}
By assumption we know that $\Ric(\omega)\leq -C^{-1} \omega$. To obtain the lower bound of $F$ we use the following argument from \cite{ruan}.
Consider a point $x$
where $F$ attains its minimum. At $x$,
$$0 \le \Delta' F \le - C^{-1} g'^{i\jb} g_{i\jb} + K_1,$$
and so $g'^{i\jb}g_{i\jb}(x) \le C$.  But the geometric-harmonic
means inequality implies that at $x$,
$$e^F = \frac{\omega'^n}{\omega^n} \ge \left( \frac{n}{g'^{i\jb}g_{i\jb}} \right)^n
\ge C,$$ and so $\sup_M |F| \le C$.  Then from Yau's estimates, as
above, we see that $\omega$ and $\omega'$ are uniformly equivalent.
Since $\Ric(\omega') \le K_2 \omega$ we see that $\de_i\de_{\jb} F \ge R_{i\jb} - K_2 g_{i \jb}$ and it follows that
$$\Delta F \le C \Delta' F + \tilde{C},$$
for a uniform constant $\tilde{C}$.
From (\ref{Delta'F}), this is uniformly bounded from above.  The
rest of the argument follows as in Lemma \ref{compactnessthm1}. Q.E.D.
\end{proof1}

\section{The Calabi flow}

In this section we prove the main theorems.
We first show that, under the assumptions of Theorems 1 and 2, the
Calabi flow exists for all time and the evolving metric  and its derivatives
are uniformly bounded.  In order to make use of the
Calabi functional, which decreases along the flow, we require the
following well-known formula \cite{calabi}, which holds for any K\"ahler metric
$\omega$:
\begin{eqnarray} \label{eqnCaformula}
\Ca(\omega) & = & \int_M \left| \Ric(\omega) - \frac{\mu}{n} \omega \right|^2
\omega^n +  \Psi,
\end{eqnarray}
where
$$\Psi = n(n-1)\pi^2\left([c_1(M)]^2 \cdot [\omega]^{n-2} - \frac{([c_1(M)]\cdot
[\omega]^{n-1})^2}{[\omega]^n} \right).$$

{\bf Remark. \cite[Corollary 1.6.3]{lazar}} Assume that $\pm
c_1(M)>0$. Then it is interesting to note that $\Psi\leq 0$ with equality if and only if
$[\omega]$ is a multiple of $c_1(M)$.  However, this inequality
has the wrong sign to be useful to us here. \\

We assume the hypotheses of Theorem 1, where $\ve$ is to be
determined later.  Let $\omega_t = \omega+ \ddbar \varphi_t$ be the
solution of the Calabi flow starting at $\omega_0$.  Applying
Lemma 1 to $\omega_0$ we immediately see that $\| \varphi_0 -
\underline{\varphi_0} \|_{C^{3+\alpha}(\omega)} \le C$ and $\omega_0
\ge C^{-1} \omega$ for any fixed $0 < \alpha < 1$.
 By the short time existence result of
\cite{che},
 the Calabi flow starting at $\omega_0$ exists on $[0,t_0]$
where $t_0>0$ is uniform, and on the subinterval $[t_0/2,t_0]$ we
have uniform $C^{k+\alpha}(\omega)$ bounds on $\varphi_t$ for all
$k\geq 4$. In particular there exist uniform constants $K_3,K_4$
such that for  $t\in [t_0/2,t_0]$ we have
\begin{equation}\label{bounds}
-K_3 \omega \leq \Ric(\omega_t)\leq K_4\omega.
\end{equation}
Now consider the set
$$\mathcal{T}:= \left\{   t\geq t_0/2\  \bigg| \ \begin{array}{l}  \textrm{the flow exists on the time interval } [0,t]  \ \textrm{and}\\
  -2K_3 \, \omega \le \Ric(\omega_{t'}) \le
2K_4 \, \omega \ \textrm{for } t' \in [t_0/2, t] \end{array}
\right\}.$$

 From what we have
just shown, there exist uniform $C^k(\omega)$ bounds along the flow for all
$k$ as long as $t$ remains in the set $\mathcal{T}$. \\

Suppose  that  $\tau=\sup \mathcal{T} < \infty$.  We will show that
there exists a uniform $\ve$ such that $\Ca(\omega_{\tau}) > \ve$,
thus contradicting assumption (ii). Since $\tau \in \mathcal{T}$, we
have uniform
 bounds for the metric $\omega_{\tau}$ and its derivatives.  At time $\tau$ there exists a point $x \in M$ at
which the following happens.
 Diagonalizing $\Ric(\omega_{\tau})$ with respect to $\omega$, one of the
 eigenvalues of $\Ric(\omega_{\tau})$ must equal $2K_4$ or $-2K_3$.  But
$\omega_{\tau}$ and $\omega$ are uniformly equivalent and  we have uniform estimates on $\omega_{\tau}$ and its derivatives, so
 there exists a uniform $\ve>0$ such that
$$\int_M |\Ric(\omega_{\tau})|^2 \omega_{\tau}^n \ge \ve,$$
and we obtain a contradiction from (\ref{eqnCaformula}) since
$\Psi=0$. \\

We turn now to the case of Theorem 2.  First, if $\delta>0$ is sufficiently
small then there exists a reference metric $\tilde{\omega}$ in $[\omega_0]$,
 uniformly equivalent to the original
$\omega$, with uniformly bounded derivatives and
$\Ric(\tilde{\omega})< - C^{-1} \tilde{\omega}$. Set $\omega_t =
\tilde{\omega} + \ddbar \varphi_t$.  From the assumptions (ii) and
(iii) we can apply Lemma 2 to obtain a uniform
$C^{3+\alpha}(\tilde{\omega})$ bound on
 $\varphi_0 -
\underline{\varphi_0}$.  We can then use a very similar argument to
the one
 above.  We replace (\ref{bounds}) by
\begin{equation}
R(\omega_{t}) \ge -K_3, \  \Ric(\omega_{t}) \le K_4 \,
\tilde{\omega}
\end{equation}
and replace $\mathcal{T}$ by
$$\mathcal{T}:= \left\{   t\geq t_0/2\  \bigg| \ \begin{array}{l}  \textrm{the flow exists on the time interval } [0,t]  \ \textrm{and}\\
 R(\omega_{t'}) \ge -2K_3, \  \Ric(\omega_{t'}) \le
2K_4 \, \tilde{\omega} \ \textrm{for } t' \in [t_0/2, t]
\end{array} \right\}.$$
If $\tau = \sup \mathcal{T} < \infty$, then at some point $x$ in $M$
either $R(\omega_{\tau}) = -2K_3$ or one of the eigenvalues of
$\Ric(\omega_t)$ with respect to $\tilde{\omega}$ is equal to
$2K_4$.  In the first alternative, since $-2K_3< \mu$, there exists a uniform $\ve$ such that $\Ca(\omega_\tau) > \ve$.  Otherwise, there exists $\ve>0$ such that
$$\int_M \left|\Ric(\omega_{\tau}) - \frac{\mu}{n} \omega_{\tau}\right|^2 \omega_{\tau}^n > 2\ve,$$
(using the fact that $\mu<0$). Picking $\delta$ small enough so that
$\Psi > - \ve$ we obtain a contradiction in the second
alternative.\\

To complete the proofs of Theorems 1 and 2, it remains to prove
convergence of the flow.  Since we have uniform $C^k$ bounds along
the flow, we know that for every sequence of times $t_k \rightarrow
\infty$ there exists a subsequence $t_{k_j}$ and a metric
$\omega_{\infty}$ such that $\omega_{t_{k_j}}$ converges in
$C^{\infty}$ to $\omega_{\infty}$ as $j\rightarrow \infty$. It is
shown in \cite{che} that such a limit $\omega_{\infty}$ must be an
extremal metric (that is, a critical metric of Calabi's functional).  But extremal metrics are unique in their K\"ahler
classes in the cases $c_1(M)=0$ \cite{calabi1954} and $c_1(M)<0$
\cite{calabi1954}, \cite{chen}. It follows that the Calabi flow
converges in $C^{\infty}$ to $\omega_{\infty}$, where $\omega_{\infty}$ is the unique  K\"ahler Ricci-flat or cscK metric in $[\omega_0]$.
  Q.E.D.

\bigskip
{\bf Acknowledgements} \ The authors would like to thank Professor S.-T. Yau, the first author's thesis advisor, for a number of very helpful discussions.  The authors are also grateful to Professor D.H. Phong for his support and advice and to Professor J. Sturm for some useful comments on an earlier version of this paper.


\end{document}